\documentclass[pdflatex,sn-mathphys]{sn-jnl}

\jyear{2021}%

\theoremstyle{thmstyleone}%
\theoremstyle{thmstyletwo}%

\theoremstyle{thmstylethree}%

\begin{document}
\title[Fourier Series in Fractional Dimensional Space]{Fourier Series in Fractional Dimensional Space}


\author*[1]{\fnm{Ali} \sur{Dorostkar}}\email{m110alidorostkar@gmail.com}

\author[1]{\fnm{Ahmad} \sur{Sabihi}}

\affil*[1]{\orgdiv{Fractal Group},  \city{Isfahan}, \country{Iran}}


\abstract{In this paper, a  Fourier series in fractional dimensional space is introduced for an arbitrarily periodic function   $f(t;\alpha)$. We call it fractional Fourier series of the order $\alpha$. Extending the basis functions of the linear space into fractional one, by rotation transformation, we define a  real and complex Fourier series and obtain their coefficients. It is also shown that the fractional derivative of a periodic function can be realized through (fractional) Fourier series with modified coefficients. }




\maketitle
\section{Introduction}\label{sec1}
Fourier series plays an important role in mathematics and applied sciences. It is a well representation of periodic function given in linear space of integer dimensions. Many people have attempted to generalize it to fractional domain. There are two methods generalizing Fourier series and discrete-time Fourier transform to fractional Fourier series and discrete-time fractional Fourier transform, respectively \cite{PS}. \\
\indent Several definitions to solve some of partial fractional differential equations \cite{HA} are developed introducing modified Riemann-Liouville fractional derivatives \cite{YC} and for solving nonlinear equations using local fractional Fourier series \cite{SH,WS,HM,YY}.  Jumarie \cite{JG1} in a paper introduces fractional sine and cosine functions and using them, he makes a Fourier series of fractional order.\\ 
However, the authors of Ref.\cite{MP} have proposed there is an invalidity since fractional Fourier series (FFS) of the basis functions $cos(n\omega t)^{\alpha}$ and $sin(n\omega t)^{\alpha}$ derived from Mittag-Leffler function might not hold. In this paper, we show how  a new definition for basis functions of fractional dimension can create the new basis functions for FFS.

\section{Preliminary and Definitions}\label{sec2}

As is well-known,  Fourier series asserts that any arbitrarily periodic function satisfying Dirichlet's sufficient conditions \cite{OAV} can be expressed by a summation of cosine and sinus functions and their coefficients are realized from the projection of function on three sets of basis functions including $1$, $cos(n\omega t)$ and $sin(n\omega t)$ (n=0,1,..), where the basis functions are orthonormal.
The Fourier series based on three sets of basis functions of $\varphi_0$, $\varphi_1$ and $\varphi_2$ can be expressed in the following where $\varphi_0$, $\varphi_1$ and $\varphi_2$ are  $1$, $cos(n\omega t)$ and $sin(n\omega t)$ (n=0,1,..), respectively. 
 \begin{eqnarray}
f(t)=\frac{A_{0}}{2}\varphi_0+\sum_{n=1}^{\infty}A_{n}\varphi_1+\sum_{n=1}^{\infty}B_{n}\varphi_2
\end{eqnarray}  
First of all, we describe a general idea for definition of the basis functions for fractional domain by the realization of linear and fractional spaces.
Based on linear algebra, a schematic of basis functions can be represented by orthogonal vectors as depicted by Fig.1 (here $\varphi_1$ and $\varphi_2$). For sake of simplicity, we describe the scenario with the two basis functions $\varphi_1$ and $\varphi_2$. Based upon a concept of linear algebra, we can move from higher (lower) dimension to lower (higher) one through the integer derivative (integral) operator. For instance, differentiating (integrating) cosine function, one can reach the sinus one and vice versa.\\
\begin{figure}
\centering
  \includegraphics[width=.5\linewidth]{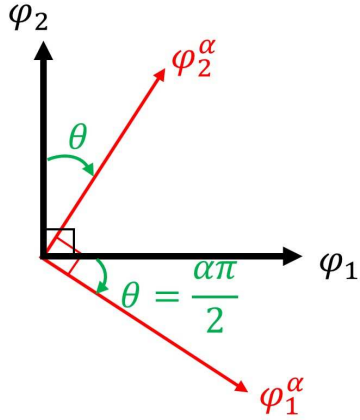}
 \caption{Schematic of basis functions for linear and fractional spaces} \label{Fig:1}
\end{figure}
\indent If suppose the fractional space is a space between two consecutive integer dimensions, then the basis functions in the fractional space can be represented as a rotation of basis functions of the linear space. Fig.1 shows that the two basis functions $\varphi_1$ and $\varphi_2$ are orthogonal to each other.  Applying the fractional derivative to $\varphi_1$ ($\varphi_2$), yields $\varphi^{\alpha}_1$ ($\varphi^{\alpha}_2$). These new basis functions are orthonormal, which make a fractional derivative to be normalized (see Definition 2).\\ 
\\
\textbf{Definition1}\\
\textit{Fractional space is a whole space between and included linear integer spaces.  }\\
\\
\textbf{Definition2}\\
\textit{Let $\varphi_0$, $\varphi_1$ and $\varphi_2$ be the basis functions of the linear space.  Let $\varphi^{\alpha}_0$,  $\varphi^{\alpha}_1$ and $\varphi^{\alpha}_2$ be the basis functions of the fractional space, then there are the following relations between the basis functions of the both spaces:
\begin{equation}
\varphi^{\alpha}_0=D^{\alpha}\bigg\{{1\bigg\}}
\end{equation}
\begin{equation}
\varphi^{\alpha}_1=\frac{1}{(n\omega)^{\alpha}}D^{\alpha}\bigg\{{cos(n\omega t)\bigg\}}=cos(n\omega t+\frac{\pi\alpha}{2})
\end{equation}
\begin{equation}
\varphi^{\alpha}_2=\frac{1}{(n\omega)^{\alpha}}D^{\alpha}\bigg\{{sin(n\omega t)\bigg\}}=sin(n\omega t+\frac{\pi\alpha}{2})
\end{equation}}\\
The factor of $\frac{1}{(n\omega)^{\alpha}}$ is used for normalization to ensure an orthonormal condition.
\\
\textbf{Lemma1}\\
\textit{  Let $\varphi^{\alpha}_0$ be a first basis function in the fractional space, then
\begin{equation}
\varphi^{\alpha}_0=D^{\alpha}\bigg\{{1\bigg\}}
=cos(\frac{\pi\alpha}{2})
\end{equation}}\\
\textbf{Proof}:\\
We consider the basis function of $\varphi_0$ in integer space based on the following expression: 
\begin{equation}
\varphi_0=\lim_{n \to 0}\bigg\{{cos(n\omega t)\bigg\}}=1
\end{equation}
Therefore, $\varphi^{\alpha}_0$ can be expressed as:
\begin{equation}
\varphi^{\alpha}_0=\lim_{n \to 0}\frac{1}{(n\omega)^{\alpha}}D^{\alpha}\bigg\{{cos(n\omega t)\bigg\}}=\lim_{n \to 0}\frac{(n\omega)^{\alpha}}{(n\omega)^{\alpha}}cos(n\omega t+\frac{\pi\alpha}{2})=cos(\frac{\pi\alpha}{2})
\end{equation}\\
\textbf{Definition3}\\
\textit{Fractional Fourier series based upon the fractional basis functions can be expressed as follows:
\begin{eqnarray}
f(t;\alpha)=a_{0}cos(\frac{\pi\alpha}{2})+\sum_{n=1}^{\infty}a_{n}cos(n\omega t+\frac{\pi\alpha}{2})+\sum_{n=1}^{\infty}b_{n}sin(n\omega t+\frac{\pi\alpha}{2})
\end{eqnarray}
where in this notation f(t; $\alpha$) denotes the function at fractional dimension $\alpha$.}\\
 Letting $\alpha=0$ changes $f(t;\alpha)$ from fractional dimension $\alpha$ into the original one in the linear space i.e. into $f(t)$ as shown in Eq.(1).
Expanding the cosine and sinus terms and rearranging, we find:\\
\begin{eqnarray}
f(t;\alpha)=a_{0}cos(\frac{\pi\alpha}{2})+\sum_{n=1}^{\infty}\bigg(a_{n}cos(\frac{\pi\alpha}{2})+b_{n}sin(\frac{\pi\alpha}{2})\bigg)cos(n\omega t)\nonumber\\+\sum_{n=1}^{\infty}\bigg(-a_{n}sin(\frac{\pi\alpha}{2})+b_{n}cos(\frac{\pi\alpha}{2})\bigg)sin(n\omega t)
\end{eqnarray}
For instance, whenever $\alpha$ is equal to $\frac{1}{2}$, the function can be expressed as follows:
\begin{equation}
f(t;\alpha =\frac{1}{2})=\frac{\sqrt2}{2}a_{0}+\frac{\sqrt2}{2}\sum_{n=1}^{\infty}(a_{n}+b_{n})cos(n\omega t)+\frac{\sqrt2}{2}\sum_{n=1}^{\infty}(-a_{n}+b_{n})sin(n\omega t)
\end{equation}

\textbf{Definition4}\\
\textit{Looking at Fig1 and the relation (9), one finds out that the FFS coefficients are a re-scaling of FS ones. This means that when moving to fractional dimension of the order $\alpha$ and projecting the function $f(t)$ on basis functions in fractional space, we really have re-scaled coefficients of FS. These re-scaled coefficients are derived from the following clockwise rotation matrix:
\begin{equation}
\left( \begin{array}{c}
A_{n}    \\
B_{n} 
\end{array} \right)=\left( \begin{array}{cc}
\cos(\frac{\pi\alpha}{2})  &   \sin(\frac{\pi\alpha}{2})  \\
-\sin(\frac{\pi\alpha}{2})  &   \cos(\frac{\pi\alpha}{2})
\end{array} \right)\left( \begin{array}{c}
a_{n}    \\
b_{n} 
\end{array} \right)
\end{equation}}\\

The FFS based upon the basis functions of the linear space can be generalized from Eq.(1) as:
 \begin{eqnarray}
f(t;\alpha)=\frac{A_{0}}{2}+\sum_{n=1}^{\infty}A_{n}cos(n\omega t)+\sum_{n=1}^{\infty}B_{n}sin(n\omega t)
\end{eqnarray}
and by a comparison of the coefficients with Eq.(9)
the  rotation matrix can be realized as follows:
\begin{equation}
R(-\frac{\pi\alpha}{2})=\left( \begin{array}{cc}
\cos(\frac{\pi\alpha}{2})  &   \sin(\frac{\pi\alpha}{2})  \\
-\sin(\frac{\pi\alpha}{2})  &   \cos(\frac{\pi\alpha}{2})
\end{array} \right)
\end{equation}\\
\textbf{Theorem1}\\
\textit{ Let the function $f(t;\alpha)$ defined on $[0,\infty)\times [0,\infty)\longrightarrow\Bbb R$ as a fractional dimension space be a continuously and integrably periodic function of the order $\alpha$ defined on $[0,\infty)\longrightarrow\Bbb R$ satisfying Dirichlet's sufficient conditions. If Fourier series of the function $f(t;\alpha)$ is expressed by the relation (8), then, the new fractional coefficients expressed by (12) are as follows:}
\begin{equation}
\left\{\begin{array}{ll}
A_{0}=2a_{0}cos(\frac{\pi\alpha}{2})  \\
A_{k}=a_{k}cos(\frac{\pi\alpha}{2})+b_{k}sin(\frac{\pi\alpha}{2} ) &    \mbox{if $k\geq 1$}\\
B_{k}=b_{k}cos(\frac{\pi\alpha}{2})-a_{k}sin(\frac{\pi\alpha}{2} )   &     \mbox{if $k\geq 1$}
\end{array} \right.
\end{equation}

\textit{and}

\begin{equation}
\left\{\begin{array}{ll}
a_{0}cos(\frac{\pi\alpha}{2})=\frac{2}{T}\int_{-\frac{T}{2}}^{\frac{T}{2}}f(t;\alpha)dt \\
a_{k}=\frac{2}{T}\int_{-\frac{T}{2}}^{\frac{T}{2}}f(t;\alpha)cos(k\omega t+\frac{\pi\alpha}{2})dt    &   \mbox{if $k\geq 1$}\\
b_{k}=\frac{2}{T}\int_{-\frac{T}{2}}^{\frac{T}{2}}f(t;\alpha)sin(k\omega t+\frac{\pi\alpha}{2})dt    &    \mbox{if $k\geq 1$}
\end{array} \right.
\end{equation}
\textbf{\textit{Proof}}\\
The new coefficients in the Fourier series expansion of $f(t;\alpha)$ given by (12) are as follows:
\begin{equation}
\begin{split}
&A_{0}=\frac{2}{T}\int_{-\frac{T}{2}}^{\frac{T}{2}}f(t;\alpha)dt=\frac{2}{T}\int_{-\frac{T}{2}}^{\frac{T}{2}}a_{0}cos(\frac{\pi\alpha}{2})dt+\frac{2}{T}\int_{-\frac{T}{2}}^{\frac{T}{2}}\sum_{n=1}^{\infty}a_{n}cos(n\omega t+\frac{\pi\alpha}{2})dt\\ &\frac{2}{T}\int_{-\frac{T}{2}}^{\frac{T}{2}}\sum_{n=1}^{\infty}b_{n}sin(n\omega t+\frac{\pi\alpha}{2})dt=2a_{0}cos(\frac{\pi\alpha}{2})
\end{split}
\end{equation}
Multiplying the both-hand sides of $f(t;\alpha)$ by $cos(k\omega t)$ and integrating over period $T$, we have
\begin{eqnarray}
A_{k}=\frac{2}{T}\int_{-\frac{T}{2}}^{\frac{T}{2}}f(t;\alpha)cos(k\omega t)dt=\frac{2}{T}\int_{-\frac{T}{2}}^{\frac{T}{2}}a_{0}cos(\frac{\pi\alpha}{2})cos(k\omega t)dt+\nonumber\\ \frac{2}{T}\int_{-\frac{T}{2}}^{\frac{T}{2}}\sum_{n=1}^{\infty}a_{n}cos(n\omega t+\frac{\pi\alpha}{2})cos(k\omega t)dt+\nonumber\\ \frac{1}{T}\int_{-\frac{T}{2}}^{\frac{T}{2}}\sum_{n=1}^{\infty}b_{n}sin(n\omega t+\frac{\pi\alpha}{2})cos(k\omega t)dt=
\end{eqnarray}  
\[ \left\{\begin{array}{ll}
0     &  \mbox{if $n\ne k$}\\
a_{k}cos\frac{\pi\alpha}{2}+b_{k}sin\frac{\pi\alpha}{2}  &   \mbox{if $n=k$}
\end{array} \right.\]
where $k\geq 1$ \\
Multiplying the both-hand sides of $f(t;\alpha)$ by $sin(k\omega t)$ and integrating over period $T$, we have 
\begin{eqnarray}
B_{k}=\frac{2}{T}\int_{-\frac{T}{2}}^{\frac{T}{2}}f(t;\alpha)sin(k\omega t)dt=\frac{2}{T}\int_{-\frac{T}{2}}^{\frac{T}{2}}a_{0}cos(\frac{\pi\alpha}{2})sin(k\omega t)dt+\nonumber\\ \frac{2}{T}\int_{-\frac{T}{2}}^{\frac{T}{2}}\sum_{n=1}^{\infty}a_{n}cos(n\omega t+\frac{\pi\alpha}{2})sin(k\omega t)dt+\nonumber\\ \frac{2}{T}\int_{-\frac{T}{2}}^{\frac{T}{2}}\sum_{n=1}^{\infty}b_{n}sin(n\omega t+\frac{\pi\alpha}{2})sin(k\omega t)dt=
\end{eqnarray}  
\[ \left\{\begin{array}{ll}
0     &  \mbox{if $n\ne k$}\\
b_{k}cos\frac{\pi\alpha}{2}-a_{k}sin\frac{\pi\alpha}{2}  &   \mbox{if $n=k$}
\end{array} \right.\]
where $k\geq 1$ \\
Note that the coefficients $a_{k}$ and $b_{k}$ are obtained from (8) when using fractional basis functions $\varphi^{\alpha}_1=cos(n\omega t+\frac{\pi\alpha}{2})$ and $\varphi^{\alpha}_2=sin(n\omega t+\frac{\pi\alpha}{2})$ of the fractional dimension space and can be realized as follows:
\begin{equation}
\left( \begin{array}{c}
a_{n}    \\
b_{n} 
\end{array} \right)=\left( \begin{array}{cc}
\cos(\frac{\pi\alpha}{2})  &   \sin(\frac{\pi\alpha}{2})  \\
-\sin(\frac{\pi\alpha}{2})  &   \cos(\frac{\pi\alpha}{2})
\end{array} \right)^{-1}\left( \begin{array}{c}
A_{n}    \\
B_{n} 
\end{array} \right)
\end{equation}\\

\section{Complex Formulation}\label{sec3}
In continue, the FFS is generalized to complex plane.

\textbf{Lemma2}\\
\textit{Let $\phi^{\alpha}_n$ be the basis function of the fractional space and $\phi_n=e^{(i n\omega t)}$ be the basis function of the linear space in the complex domain, then 
\begin{equation}
\phi^{\alpha}_n=\frac{1}{(n\omega)^{\alpha}}D^{\alpha}\bigg\{\phi_n\bigg\}=e^{i(n\omega t+\frac{\pi\alpha}{2})}
\end{equation}}
 
\textbf{\textit{Proof}}\\
The proof is easily made as follows.
 $\phi^{\alpha}_n$ for $n> 0$ and $n<0$ can be calculated as:
\begin{equation}
\begin{split}
\phi^{\alpha}_n=\frac{1}{(n\omega)^{\alpha}}D^{\alpha}\bigg\{\phi_n\bigg\}=\frac{1}{(n\omega)^{\alpha}}D^{\alpha}\bigg\{e^{(in\omega t)}\bigg\}=\frac{1}{(n\omega)^{\alpha}}(in\omega)^{\alpha}e^{(in\omega t)}\\\frac{1}{(n\omega)^{\alpha}}(n\omega)^{\alpha}(i^{\alpha})e^{(in\omega t)}=(i^{\alpha})e^{(in\omega t)}=(e^{(\frac{i\pi}{2})})^{\alpha}e^{(in\omega t)}=e^{i(n\omega t+\frac{\pi \alpha}{2})}
\end{split}
\end{equation}
and for $\phi^{\alpha}_0$, we have:
\begin{eqnarray}
\phi^{\alpha}_0=D^{\alpha}\bigg\{1\bigg\}=D^{\alpha}\bigg\{{\lim_{n \to 0}}e^{+in\omega t}\bigg\}=e^{+i\frac{\pi\alpha}{2}}
\end{eqnarray} 
Based on Eqs. (3), (4) and (5) the basis functions $\varphi^{\alpha}_i (i=0,1,2)$ and $\phi^{\alpha}_n$  in real and complex domains respectively can be related to each other as follows:
\begin{eqnarray}
\varphi^{\alpha}_0=\frac{\phi^{\alpha}_{0}+(\phi^{\alpha}_{0})^{*}}{2}=\Re\{\phi^{\alpha}_{0}\} 
\end{eqnarray} 
\begin{eqnarray}
\varphi^{\alpha}_1=\frac{\phi^{\alpha}_n+(\phi^{\alpha}_n)^{*}}{2}=\Re\{\phi^{\alpha}_{n}\} 
\end{eqnarray} 
\begin{eqnarray}
\varphi^{\alpha}_2=\frac{\phi^{\alpha}_n-(\phi^{\alpha}_n)^{*}}{2i}=\Im\{\phi^{\alpha}_{n}\}  \end{eqnarray} 
where $\Re$ and $\Im$ are real and imaginary part, respectively.
Therefore the general form of basis function can be realized as:
\begin{eqnarray}
\phi^{\alpha}_n=e^{i(n\omega t+\frac{\pi \alpha}{2})}
\end{eqnarray}
It can be observed that in-spite of validity properties of $\varphi_{-n}=(\varphi_{n})^{*}$  in conventional FS  ,however;  this type of properties in fractional space is invalid.  
\begin{eqnarray}
\phi^{\alpha}_{-n}=e^{i(-n\omega t+\frac{\pi \alpha}{2})}\neq (\phi^{\alpha}_{n})^*
\end{eqnarray}
\\
\textbf{Theorem2}\\
\textit{Let the complex function $f(t;\alpha)=\sum_{n=-\infty}^{+\infty}c^{\alpha}_{n}e^{i(n\omega t+\frac{\pi\alpha}{2})}$ satisfy Dirichlet's sufficient conditions and be defined on $[0,\infty)\times [0,\infty)\longrightarrow\Bbb C$, then implies that\\
\begin{eqnarray}
f(t;\alpha)=\sum_{n=-\infty}^{+\infty}c^{\alpha}_{n}e^{i(n\omega t+\frac{\pi\alpha}{2})}=a_{0}cos(\frac{\pi\alpha}{2})+\sum_{n=1}^{\infty}a_{n}cos(n\omega t+\frac{\pi\alpha}{2})+\nonumber\\ \sum_{n=1}^{\infty}b_{n}sin(n\omega t+\frac{\pi\alpha}{2})~~~~~~~~~~~~~~~~~~~~~~~~~
\end{eqnarray}
where the complex coefficients are:}\\
\begin{equation}
\left\{\begin{array}{ll}
c_{0}^{\alpha}=a_{0}cos(\frac{\pi\alpha}{2})e^{-\frac{i\pi\alpha}{2}}  &  \mbox{ $n=0$}\\
c_{n}^{\alpha}=\frac{1}{2}(a_{n}-ib_{n})  &   \mbox{if $n\geq 1$}\\
c_{-n}^{\alpha}=\frac{1}{2}(a_{-n}+ib_{-n})e^{-i\pi\alpha}    &    \mbox{if $n\leq -1$}
\end{array} \right.
\end{equation}
\textbf{\textit{Proof}}\\
The proof is easily made as follows:\\
For all $n\in \Bbb Z$, we expand the right-hand side of (28) and show that it is identical to the left one. Expanding the right-hand side gives us 
\begin{eqnarray}
f(t;\alpha)=a_{0}cos(\frac{\pi\alpha}{2})+\sum_{n=1}^{\infty}a_{n}cos(n\omega t+\frac{\pi\alpha}{2})+\sum_{n=1}^{\infty}b_{n}sin(n\omega t+\frac{\pi\alpha}{2})=\nonumber\\
a_{0}cos(\frac{\pi\alpha}{2})+\sum_{n=1}^{\infty}a_{n}\bigg\{\frac{e^{i(n\omega t+\frac{\pi\alpha}{2})}+e^{-i(n\omega t+\frac{\pi\alpha}{2})}}{2}\bigg\}+~~~~~~~~~~~~~~~~~~~\nonumber\\
\sum_{n=1}^{\infty}b_{n}\bigg\{\frac{e^{i(n\omega t+\frac{\pi\alpha}{2})}-e^{-i(n\omega t+\frac{\pi\alpha}{2})}}{2i}\bigg\}=~~~~~~~~~~~~~~~~~~~~\nonumber\\ a_{0}cos(\frac{\pi\alpha}{2})+\sum_{n=1}^{\infty}(\frac{a_{n}}{2}+\frac{b_{n}}{2i})e^{i(n\omega t+\frac{\pi\alpha}{2})}+\sum_{n=1}^{\infty}(\frac{a_{n}}{2}-\frac{b_{n}}{2i})e^{-i(n\omega t+\frac{\pi\alpha}{2})}~~~~~~~
\end{eqnarray}

Changing variable $n$ by $-N$ at the right term of the right-hand side of (30), we have
\begin{eqnarray}
f(t;\alpha)=a_{0}cos(\frac{\pi\alpha}{2})+\sum_{n=1}^{\infty}(\frac{a_{n}}{2}+\frac{b_{n}}{2i})e^{i(n\omega t+\frac{\pi\alpha}{2})}+~~~~~~~~~~~~~~~\nonumber\\\sum_{N=-\infty}^{N=-1}(\frac{a_{-N}}{2}-\frac{b_{-N}}{2i})e^{-i(-N\omega t+\frac{\pi\alpha}{2})}=~~~~~~~~~~~~~~~~~\nonumber\\ a_{0}cos(\frac{\pi\alpha}{2})+\sum_{n=1}^{\infty}(\frac{a_{n}}{2}+\frac{b_{n}}{2i})e^{i(n\omega t+\frac{\pi\alpha}{2})}+\sum_{n=-\infty}^{n=-1}(\frac{a_{-n}}{2}-\frac{b_{-n}}{2i})e^{-i\pi \alpha}e^{i(n\omega t+\frac{\pi\alpha}{2})}=\nonumber\\
\sum_{n=-\infty}^{+\infty}c^{\alpha}_{n}e^{i(n\omega t+\frac{\pi\alpha}{2})}~~~~~~~~~~~~~~~~~~~~~~~~~~~~~~~~~~~
\end{eqnarray}
Getting $a_{0}cos(\frac{\pi\alpha}{2})=c_{0}^{\alpha}e^{\frac{i\pi\alpha}{2}}$ for $n=0$, the right-hand side of (31) would imply the left hand-side of (28).\\ 
Comparing the both sides of (31) together, we find the coefficients given by (29).

\section {The relations between Fractional Derivative, FS and FFS}
If we take a fractional derivative from Eq.(1) , then we have
  \begin{eqnarray}
D^{\alpha}f(t)=\frac{A_{0}}{2}D^{\alpha}\big\{\varphi_0 \big\}+\sum_{n=1}^{\infty}A_{n}D^{\alpha}\big\{\varphi_1 \big\}+\sum_{n=1}^{\infty}B_{n}D^{\alpha}\big\{\varphi_2 \big\}
\end{eqnarray} 
 and consequently by applying the direct fractional derivative of basis functions of cosine and sinus functions, we have:
 \begin{eqnarray}
D^{\alpha}f(t)= f(t;\alpha)=\frac{A_{0}}{2}D^{\alpha}\big\{1 \big\}+\sum_{n=1}^{\infty}A_{n}cos(n\omega t+\frac{\pi} {2}\alpha)+\sum_{n=1}^{\infty}B_{n}sin(n\omega t+\frac{\pi}{2}\alpha)\nonumber\\
=\frac{A_{0}}{2} cos(\frac{\pi} {2}\alpha)+\sum_{n=1}^{\infty}A_{n}cos(n\omega t+\frac{\pi} {2}\alpha)+\sum_{n=1}^{\infty}B_{n}sin(n\omega t+\frac{\pi}{2}\alpha)~~~~~~~~~~~~~~~~~~~~
\end{eqnarray}
then, based on the theorem 1 it can be figured out in terms of a conventional FS with new modified coefficients.
 \begin{eqnarray}
D^{\alpha}f(t)= f(t;\alpha)=\frac{A'_{0}}{2} +\sum_{n=1}^{\infty}A'_{n}cos(n\omega t)+\sum_{n=1}^{\infty}B'_{n}sin(n\omega t)\nonumber\\
\end{eqnarray}
where 
\begin{equation}
\left\{\begin{array}{ll}
A'_{0}=A_{0}cos(\frac{\pi\alpha}{2})  \\
A'_{k}=A_{n}cos(\frac{\pi\alpha}{2})+B_{n}sin(\frac{\pi\alpha}{2} ) &    \mbox{if $n\geq 1$}\\
B'_{k}=B_{n}cos(\frac{\pi\alpha}{2})-A_{n}sin(\frac{\pi\alpha}{2} )   &     \mbox{if $n\geq 1$}
\end{array} \right.
\end{equation}



\begin{thebibliography}{99}

\bibitem{HA} M.A. Hammad and R.R. Khalil,"Fractional Fourier Series with Applications", Amer.J.Comput. Appl. Math., 4 (6), 187-191, 2014.

\bibitem{HM} M.-S. Hu, R. P. Agarwal, and X.-J. Yang, “Local fractional Fourier series with Application to Wave Equation in Fractal Vibrating String,” Abs. Appl. Analysis, 2012, Article ID 567401, 15 pages, 2012.

\bibitem{JG1} G. Jumarie, "Fourier’s transform of fractional order via Mittag-Leffler function and modified Riemann-Liouville derivative," J. Appl. Math. Informatics, 26 (5–6), 1101-1121, 2008.

\bibitem {MP}  P. R. Massopust and Ahmed I. Zayed. "On the Invalidity of Fourier series Expansions of Fractional Order." Fract. Calc. and Appl. Analysis, 18 (6),1507-1517, 2015.

\bibitem {OAV} A.V. Oppenheim, and A.S. Willsky, Signals and Systems, Prentice Hall, p.198, 1997 

\bibitem{PS} S.C. Pei, M.H. Yeh, and T.L. Luo, "Fractional Fourier Series Expansion for
Finite Signals and Dual Extension to Discrete-Time Fractional Fourier Transform", IEEE Trans.  Sig. Process. 47 (10), 2883-2888, 1999.

\bibitem{SH} H. M. Srivastava, A. K. Golmankhaneh, D. Baleanu, and X.-J.Yang, “Local Fractional Sumudu Transform with Application to IVPs on Cantor Sets,” Abs. and Appl. Analysis, vol. 2014, Article ID 620529, 7 pages, 2014.

\bibitem{WS} S.-Q. Wang, Y.-J. Yang, and H. K. Jassim, “Local Fractional
Function Decomposition Method for Solving Inhomogeneous Wave Equations with Local Fractional Derivative,” Abs. and Appl. Analysis, 2014, Article ID 176395, 7 pages, 2014.

\bibitem{YY} Y.J. Yang and S.Q. Wang "Local Fractional Fourier Series Method for Solving Nonlinear Equations with Local Fractional Operators", Math. Prob. Eng. 2015, Article ID: 481905, 9 pages, 2015. 

\bibitem{YC} C.Yu, "Fractional Fourier Series and its Applications", J. Physics: conference series, 1976, 2021.


\end{thebibliography}

\end{document}